\documentclass[11pt]{article}
\usepackage{latexsym}

\title{Compact endomorphisms of Banach algebras of infinitely differentiable
functions}
\author{Joel F. Feinstein and Herbert Kamowitz}

\begin{document}
\newcommand{\dis}{\displaystyle}
\newcommand{\pis}{\mbox{$\overline{\phi(b)}$}}

\newcommand{\QED}{\hfill$\Box$\vskip.3cm}
\def\unorm#1#2{{\|#1\|_{{#2}}}}
\def\norm#1{{\|#1\|}}
\def\rd{{\rm d}}

\maketitle

\begin{abstract}
Let $(M_n)$ be a sequence of positive numbers
satisfying $M_0=1$ and 
$\dis \frac{M_{n+m}}{M_n M_m} \geq{{n+m}\choose{m}}$
for all non-negative integers $m$, $n$.
We let
\[D([0,1], M)=\{ f \in C^{\infty}([0,1]):\|f\|_{D}=\sum_{n=0}^{\infty}
\frac{\|f^{(n)}\|_{\infty}}{M_n} < \infty\}.\]
    With pointwise addition and multiplication, $D([0,1],M)$ is a
unital commutative semisimple Banach algebra.  If $\lim_{n\to\infty}
(n!/M_n)^{1/n}=0,$ then the maximal ideal space of the algebra
is $[0,1]$ and every non-zero endomorphism $T$ has the form
$Tf(x)=f(\phi(x))$ for some selfmap $\phi$ of the unit interval.
Previously we have shown for a wide class of $\phi$ mapping the unit interval to
itself that if $\|\phi'\|_\infty < 1$ then $\phi$ induces a compact
endomorphism. Here we investigate the extent to which this condition is necessary, and
we determine the spectra of all compact endomorphisms of $D([0,1],M).$
We also simplify and strengthen some of our earlier
results on general endomorphisms of $D([0,1],M)$.
\end{abstract}

{\bf Introduction}

    In two previous  papers  \cite{jfhk} and  \cite{blaub}
we studied endomorphisms
of certain algebras of infinitely differentiable functions.
In this
note, we continue our study, concentrating on algebras of functions
on compact intervals in the real line.
Let $(M_n)$ be a sequence of positive numbers
satisfying $M_0=1$ and $\dis \frac{M_{n+m}}{M_nM_m} \geq
\left(\begin{array}{c}n+m\\n \end{array}\right)$, $m$, $n$,
non-negative integers. If $X$ is a compact interval of positive length
we let
\[D(X, M)=\{ f \in C^{\infty}(X):\|f\|_{D}=\sum_{n=0}^{\infty}
\frac{\|f^{(n)}\|_{\infty}}{M_n} < \infty\}.\]

    With pointwise addition and multiplication, $D(X,M)$ is a
unital commutative semisimple Banach algebra.
The maximal ideal
space of $D(X,M)$ depends on $(M_n).$
In particular, it was shown in \cite{dd}
(see also \cite{dales} Theorem 4.4.16)
that if
$(M_n)$ is a {\it non-analytic} sequence, i.e. $\dis \lim_{n \rightarrow
\infty}
(\frac{n!}{M_n})^{1/n}=0$, then the maximal ideal space of
$D(X,M)$ is precisely $X$.

A stronger condition on the sequence $(M_n)$ is that
$\dis \sum \frac{M_n}{M_{n+1}} < \infty.$
In this case the sequence $(M_n)$ is non-analytic and, in addition,
$D(X,M)$ is regular (\cite{dales}, Theorem 4.4.22).
Such a weight sequence will be called
{\it non-quasi-analytic.}

When the weight
$(M_n)$ is non-analytic, every non-zero endomorphism $T$ of the
algebra $D(X,M)$ has the form $Tf(x)=f(\phi(x))$ for
some $\phi:X \rightarrow X.$
    This leads us to ask which $\phi:X \rightarrow X$ induce
endomorphisms of $D(X,M)$?
 That is, which $\phi$ satisfy $f \circ \phi \in D(X,M)$
whenever $f \in D(X,M)$?
 We observe that since $f(x)=x$ is in
$D(X,M)$ for all $X$ and $(M_n)$, for $\phi$ to induce an endomorphism
it is necessary that $\phi \in D(X,M).$

We sometimes require an additional condition on an infinitely differentiable selfmap $\phi$
of the interval $X$: we say that
$\phi$ is {\it analytic} if
$$\dis \sup_{k}\left({\frac{\|\phi^{(k)}\|_\infty}{k!}}\right)^{1/k} < \infty.$$

In \cite{blaub} the following theorems were proved concerning analytic
    selfmaps $\phi$ of $[0,1]$.

{\bf Theorem A:} Let $(M_n)$ be a non-analytic weight sequence and let $\phi$ be an analytic
selfmap of $[0,1]$
with $\|\phi'\|_\infty < 1.$
Then $\phi$ induces an endomorphism of $D([0,1],M)$, and this endomorphism is compact.

{\bf Theorem A$'$:} Let $(M_n)$ be a non-analytic weight sequence such that there is a constant
$B>0$ with
$\dis \frac{M_m}{m!}\frac{n!}{M_n} \leq \frac{B}{m^{n-m}}$. Let $\phi$ be an infinitely differentiable
selfmap of $[0,1]$ such that $\|\phi'\|_\infty \leq 1$,
and $\dis (\frac{\|\phi^{(k)}\|_\infty}{k!})$ is bounded. Then $\phi$ induces an
endomorphism of  $D([0,1],M).$

{\bf Theorem B:} Suppose that $(M_n)$ is a non-analytic weight sequence and that
$\phi$ is an analytic selfmap of $[0,1]$. If $\phi$
induces an endomorphism of $D([0,1],M)$, then $\|\phi'\|_\infty \leq 1.$

In Part I of this paper we investigate when a converse to the compactness part of
Theorem A is true for
$D([0,1],M);$ namely, if $\phi$ induces a compact endomorphism of
$D([0,1],M)$ must $\|\phi'\|_\infty < 1$?
Along the way we prove a stronger version
of Theorem B, eliminating the analyticity assumption on $\phi$
and removing the need for the technical calculations used in the original proof.
In Part II, we determine the spectra of compact endomorphisms of $D([0,1],M).$
Part III contains results on general endomorphisms of $D([0,1],M).$ In particular,
we simplify and weaken
slightly the growth condition on the sequence $M$ from Theorem A$'$ and we investigate
the extent to which this condition is necessary.

We would like to thank the referee for some very helpful comments.

{\bf Part I. Compact endomorphisms}

    We begin by recalling some facts about endomorphisms.
Let $B$ be a commutative, unital semisimple Banach algebra
with a connected maximal ideal space $X$,
and suppose that
$T$ is a non-zero endomorphism of $B$. Then there exists a
$w*$-continuous function $\phi:X \rightarrow X$
with $\widehat{Tf}(x)=\hat{f}(\phi(x))$
for $f \in B$ and $x \in X.$
We let $\phi_n$ denote the $n^{\rm th}$ iterate of $\phi$.
If the endomorphism $T$ is compact, then it was
shown in \cite{comp1} that
$\dis \bigcap_{n=0}^{\infty}{ \phi_n(X)}=\{x_0\}$ where
$x_0 \in X$ is a fixed point of $\phi$ which clearly must be unique.
In view of this the following result is not surprising.

{\bf Theorem 1.1:} Suppose that $X$ is a  compact interval of positive length
and $(M_n)$ is a weight sequence. Suppose that
$\phi$ induces a compact endomorphism
of $D(X,M)$. If $x_0$ is the fixed point of $\phi$, then
$|\phi'(x_0)|<1.$

{\bf Proof:} Suppose, for a contradiction, that $|\phi'(x_0)| \geq 1.$
Then $\dis \|\phi_n\|_{D} \geq
\frac{|\phi_n'(x_0)|}{M_1}=\frac{(|\phi'(x_0)|)^n}{M_1} \geq \frac{1}{M_1}.$
It follows that the sequence of iterates $\phi_n$ has a subsequence
$\phi_{n_k}$ for which there is a $\delta > 0$
with $\dis \frac{\|\phi_{n_k+1}\|_{D}}{\|\phi_{n_k}\|_{D}}
\geq \delta.$
   For this subsequence let $\dis f_k=
\frac{\phi_{n_k}}{\|\phi_{n_k}\|_{D}}.$ Since $\|f_k\|_{D}=1$ 
and $\phi$
induces a compact endomorphism there exists $g \in D(X,M)$
and a subsequence $\dis \{f_{k_j}\}$ with $\dis f_{k_j}(\phi)
\rightarrow g$ in norm. With no loss of generality we call the
subsequence $\dis \{f_k\}.$ Then $\dis f_k(\phi)=\frac{\phi_{n_k+1}}
{\|\phi_{n_k}\|_{D}} \rightarrow g$ in norm.
Since $\dis \bigcap^\infty_{n=0} \phi_n(X)=\{x_0\},$ we have
 $\dis \lim_{n \rightarrow \infty} \phi_n(x)=x_0$
for each $x \in X.$ Now
\[g(x_0)=\lim_{k\rightarrow \infty} f_k(\phi(x_0))=\lim_{k \rightarrow
\infty}
\frac{\phi_{n_k+1}(x_0)}{\|\phi_{n_k}\|_{D}}=\lim_{k \rightarrow \infty}
\frac{x_0}{\|\phi_{n_k}\|_{D}}.\]
Suppose, at first, that $x_0 \neq 0.$ Then $\dis \lim_{k \rightarrow \infty}
\frac{1}{\|\phi_{n_k}\|_{D}}=\frac{g(x_0)}{x_0}.$ 
Hence for each $x \in X,$
 \[g(x)=\lim_{k \rightarrow \infty}f_k(\phi(x))=\lim_
{k \rightarrow \infty} \frac{\phi_{n_k+1}(x)}{\|\phi_{n_k}\|_{D}}
=g(x_0),\]
whence $g$ is a constant function.

     We consider in turn the two cases $g(x_0)=0$ and $g(x_0) \neq 0.$

(a) First, if $g(x_0)=0,$ then $\dis \lim_{k \rightarrow \infty}
f_k(\phi)=0$
in norm. But $\dis
\|f_k(\phi)\|_{D}=\frac{\|\phi_{n_k+1}\|_{D}}{\|\phi_{n_k}\|_{D}}
\geq \delta > 0$, a contradiction.

(b) Otherwise, $g(x_0) \neq 0.$ Since $\dis f_k(\phi) \rightarrow g$
in norm, we have $(f_k \circ \phi)'(x) \rightarrow g'(x)$ for all
$x \in X.$ Since $g$ is a constant function, $g'(x)=0$ for all
$x \in X.$ However,
\[\lim_{k \rightarrow \infty}|(f_k \circ \phi)'(x_0)|=
\lim_{k \rightarrow \infty}|\frac{\phi_{n_k+1}'(x_0)}{\|\phi_{n_k}\|_{D}}|
\geq
\lim_{k \rightarrow \infty} \frac{1}{\|\phi_{n_k}\|_{D}}=\frac{g(x_0)}{x_0}
\neq 0,\]
again a contradiction. Hence, if $|\phi'(x_0)| \geq 1$ and $x_0 \neq 0,$ then
$\phi$ does not
induce a compact endomorphism.


  Finally, if $x_0=0$, we have
 \[g(x)=\lim_{k \rightarrow \infty}f_k(\phi(x))=\lim_
{k \rightarrow \infty} \frac{\phi_{n_k+1}(x)}{\|\phi_{n_k}\|_{D}}
=0\]
for all $x$, and  as in (a) above we again are led to a contradiction. Thus if
$|\phi'(x_0)| \geq 1,$ then
 $\phi$ does not induce a compact endomorphism.
\QED

   We remark that this result shows that if $\phi$ induces a compact
endomorphism of $D([0,1],M)$, then some iterate $\phi_n$ satisfies
$\|\phi_n'\|_\infty < 1.$ We conjecture that
$\phi$ itself satisfies $\|\phi'\|_\infty < 1.$ Most of the remainder of this
section is devoted to this question. Before we get on to the main focus
of this section, however, we use Theorem 1.1 to produce a very short proof of
a stronger version of Theorem B. In particular, observe that there is no
restriction on the growth of the derivatives of $\phi$.

{\bf Theorem 1.2:} Let $(M_n)$ be a non-analytic weight and let $\phi$
be a selfmap of $[0,1].$
If $\phi$
induces an endomorphism of $D([0,1],M)$, then $\|\phi'\|_\infty \leq 1.$

{\bf Proof:} Let $\phi$ induce an endomorphism of $D([0,1],M)$ and
suppose that $|\phi'(x_0)| > 1$ for some $x_0 \in [0,1].$ By continuity
of $\phi'$ it is no restriction to assume that $0 < x_0 <1.$
It is easy to see that there exists a polynomial function $p:[0,1]
\rightarrow
[0,1]$ such that $p(\phi(x_0)) = x_0$,
$\dis p'(\phi(x_0))=\frac{1}{|\phi'(x_0)|}$, and $\|p'\|_\infty < 1.$
It then follows from Theorem A that $p$ induces a
compact endomorphism of $D([0,1],M)$, whence $p \circ \phi$ also induces
a compact endomorphism of the algebra. But $(p \circ \phi)(x_0)=x_0$,
and  $|(p \circ \phi)'(x_0)|=|p'(\phi(x_0))\phi'(x_0)|=1.$ Since
$x_0$ is a fixed point of $p \circ \phi$ we obtain a contradiction to
Theorem 1.1. Thus $\|\phi'\|_\infty \leq 1.$
\QED


  We now continue with some results about homomorphisms
rather than endomorphisms.


Let $X$, $Y$ be compact intervals in $\bf R$ of positive length, let $\phi$
be a map from $X$ to $Y$ and let
$(M_n)$ be a non-analytic weight sequence. For $f$ in $D(Y,M)$
we  look at $f \circ \phi$
and investigate whether it is in $D(X,M)$.
 In other words we  ask whether $\phi$
induces a homomorphism from $D(Y,M)$ to $D(X,M)$. (Since these algebras are
commutative,
semi-simple Banach algebras,
if this happens then the homomorphism must be continuous.)
We also ask whether the homomorphism is compact.

We begin with some elementary observations, which follow easily from
properties of the relevant restriction maps.


{\bf Observation 1} Suppose that
$\phi:X \rightarrow Y$ induces a homomorphism from $D(Y,M)$ to $D(X,M)$.
Let $X_0$ be a closed interval of positive length contained in $X.$
Then $\phi_{|X_0}$ induces a homomorphism from $D(Y,M)$ to $D(X_0,M).$
Moreover, if the homomorphism induced by $\phi$ is compact, then
so is the homomorphism induced by $\phi_{|X_0}.$

{\bf Observation 2:}
Let $Y_0$ be a closed interval in $Y$ such that $\phi(X) \subseteq Y_0$.
If $\phi$ induces a homomorphism from $D(Y_0,M)$ to $D(X,M)$ then
$\phi$ induces a homomorphism from $D(Y,M)$ to $D(X,M)$.
Similarly,
If $\phi$ induces a compact homomorphism from $D(Y_0,M)$ to $D(X,M)$ then
$\phi$ induces a compact homomorphism from $D(Y,M)$ to $D(X,M)$.

We now
prove the converse to the second of these observations
under the additional assumption that the algebra is regular.
We do not know whether this assumption is redundant.

{\bf Lemma 1.3} If $M$ is a non-quasi-analytic sequence (so that $D(Y,M)$ is
regular) and $\phi(X)$ is a subset of the interior (relative to $Y$) of $Y_0$
then
$\phi$ induces a homomorphism from $D(Y_0,M)$ to $D(X,M)$ if and only if
$\phi$ induces a homomorphism from $D(Y,M)$ to $D(X,M)$, and the compactness
of either homomorphism implies that of the other.

{\bf Proof:}
Suppose that $D(Y,M)$ is regular and that $\phi(X)$
is a subset of the interior (relative to $Y$) of $Y_0$.
By regularity we may choose functions $g$, $h$ in $D(Y,M)$ such that
$g$ is supported on the interior of $Y_0$, the support of $h$ does not meet
$\phi(X)$ and $g+h = 1$. Write $g_1=g_{|Y_0}$ and $h_1=h_{|Y_0}$.
Suppose that $\phi$ induces a homomorphism from $D(Y,M)$ to $D(X,M)$.
Let $f \in D(Y_0,M)$. We have $f=fg_1 + fh_1$ and so (since
$h_1 \circ \phi = 0$)
$f \circ \phi = (fg_1) \circ \phi$. Since $fg_1$ is supported
on the interior of $Y_0$ we may extend it to obtain a function $F$
in $D(Y,M)$ with $F\circ \phi = (fg_1) \circ \phi = f \circ \phi $.
By the assumption on $\phi$, $F \circ \phi \in D(X,M)$, i.e.
$f \circ \phi \in D(X,M)$. This shows that
$\phi$ induces a homomorphism from $D(Y_0,M)$ to $D(X,M)$.
Next suppose that the homomorphism induced from $D(Y,M)$ to $D(X,M)$ is
compact. Let $f_n$ be a bounded sequence of functions in $D(Y_0,M)$.
As above we have $f_n = f_n g_1 + f_n h_1$ and
$f_n \circ \phi = (f_n g_1) \circ \phi$.
Since the supports of the bounded sequence of functions $f_n g_1$
are contained in the interior of $Y_0$, these functions may be extended
to functions $F_n$ in $D(Y,M)$ with the same norm and such that
$F_n \circ \phi = (f_n g_1)\circ \phi = f_n \circ \phi$.
Now, by the compactness assumption, $F_n \circ \phi$ has a convergent
subsequence in $D(X,M)$, i.e. $f_n \circ \phi$ has a convergent subsequence.
This shows that the homomorphism from $D(Y_0,M)$ to $D(X,M)$ induced by
$\phi$ is compact.

Combining this with the second observation above gives the result.
\QED

The next lemma may appear to be obvious. Note, however,
that our proof is only valid in the
case where the algebra is regular.

{\bf Lemma 1.4:} Suppose that $(M_n)$ is non-quasi-analytic and that
$\phi'$ is constantly $1$ on $X$.
Then $\phi$ does not induce a compact homomorphism from $D(Y,M)$ to
$D(X,M)$.
\par
{\bf Proof:}
We may assume (after translating)
that $\phi(x)=x$. Choose a non-degenerate compact interval $X_0$ contained in
the interior of $X$ and let $E_1$, $E_2$ be the closed subspaces of,
respectively, $D(X,M)$, $D(Y,M)$ consisting of those elements whose supports
are subsets of $X_0$. Then these infinite-dimensional Banach spaces are
isometrically isomorphic to each other using the restriction to $E_2$
of the homomorphism induced by $\phi$. Since this map is clearly not
compact, neither is the original homomorphism.
\QED

We are now ready to give our main result concerning necessary
conditions on
$\phi$ in order for $\phi$ to induce a compact endomorphism.
Note that we do not assume that $\phi$ is analytic here.

{\bf Theorem 1.5:}
Let $(M_n)$ be a non-quasi-analytic algebra sequence and let
$\phi$ be a map from $[0,1]$ to $[0,1]$. Suppose that $\phi$ induces a
compact endomorphism of $D([0,1],M)$. Then for all $x \in (0,1)$ we have
$|\phi'(x)| < 1$. Thus, if $\phi$ induces a compact endomorphism of
$D([0,1],M)$ and $|\phi'|$ peaks in $(0,1),$ then $\|\phi'\|_\infty < 1.$


{\bf Proof:}
  Suppose that $\phi$ induces a compact endomorphism of $D([0,1],M)$.
Since $\phi$ induces an endomorphism
of $D([0,1],M)$ it follows from Theorem 1.2 that $\|\phi'\|_\infty \leq 1.$
Now suppose, for contradiction, that there is an $x_1 \in (0,1)$ with
$|\phi'(x_1)|=1$. Composing with the map $t \mapsto (1-t)$ if necessary,
we may assume that $\phi'(x_1) = 1$. Certainly
$\phi(x_1) \in (0,1)$. Set
$$\delta = \min\{x_1, 1-x_1, \phi(x_1), 1-\phi(x_1)\}.$$
Set $a=x_1-\delta$, $b=x_1+\delta$, $c=x_1-\phi(x_1)$,
$X=[a,b]$ and $Y = [c, c+1]$.
Define $\psi:X \rightarrow Y$ by
$\psi(x)=\phi(x)+ c$. From
our assumptions, the observations above and obvious properties of
translations,
it is evident that
$\psi$ must induce a compact homomorphism from $D(Y,M)$ to $D(X,M)$.

Set $Y_0=X$. Then $\psi(X) \subseteq Y_0 \subseteq Y$.
Note that $\psi(x_1)=x_1$ and $\psi'(x_1)=1$ so, by Theorem 1.1,
$\psi$ does not induce a compact endomorphism of $D(X,M)$ i.e.
$\psi$ does not induce a compact homomorphism from $D(Y_0,M)$ to
$D(X,M)$. In view of the observations above,
this will contradict Lemma 1.3 provided that $\psi(X)$ is a subset of the
interior of $Y_0$. The only way this could fail would be if $\psi'$
was constantly equal to $1$ on a non-degenerate compact
subinterval of $X$, say
$X_0$.
This, however, would contradict the observations above,
since, by Lemma 1.4, $\psi|{X_0}$
would not induce a compact homomorphism
from $D(Y,M)$ to $D(X_0,M)$. The
result now follows.
\QED

 The only remaining case is when $\phi'(0)=1$ or $\phi'(1)=1$.
Theorem 1.1 deals with this if the relevant
point is a fixed point.
Using the flip of the interval, the only remaining case is really
when $\phi'(1) = 1$ and $\phi(1) \in (0,1)$.
That is, can we find a weight sequence $(M_n)$ and an analytic $\phi \in
D([0,1],M)$ with $\phi'(1)=1$ and $0 < \phi(1) < 1$ such that $\phi$ induces
a compact endomorphism of $D([0,1],M)$?
This question is still open.

   However, we have two partial results whose proofs we omit.
   First we can show that if $\phi$ is a quadratic function, and if
$\phi$ induces a compact endomorphism of $D([0,1],M)$,
then $\|\phi'\|_\infty < 1.$ Secondly, the function $g$ defined by
$\dis g(z)=\sum_{n=0}^\infty \frac{z^n}{M_n}$ is an
entire function. Using more sophisticated results on entire functions
as developed
in \cite{boas} and \cite{La}, we are also able to show that
if the order of
$g$ is less than $\frac{1}{2}$, or if the order of $g$ equals
$\frac{1}{2}$ and the type of $g$ is finite,
and if an analytic selfmap
$\phi$ of $[0,1]$ induces a compact endomorphism of $D([0,1],M),$ then again
$\|\phi'\|_\infty < 1.$

{\bf Part II. Spectra of compact endomorphisms}

    We next turn to the question of determining the spectra
of compact endomorphisms of $D([0,1],M)$.
    We denote the spectrum of an operator $T$ by $\sigma(T)$. We
also note that if $T$ is a non-zero endomorphism of $D([0,1],M)$, then
$\lambda=1 \in \sigma(T)$ since $T1=1.$ Also, if $T$ is a compact
endomorphism, then $0 \in \sigma(T).$

We will show (in Theorems 2.4 and 2.5)
that if $\phi$ induces a compact endomorphism $T$ of
$D([0,1],M)$ and if $x_0$ is the fixed point of $\phi$, then
$\dis \sigma(T)=\{(\phi'(x_0))^n: n$ is a positive integer$\} \cup \{0,1\},$
and each non-zero element of $\sigma(T)$ is an eigenvalue of multiplicity $1.$
First we require some lemmas.

{\bf Lemma 2.1:} Let $(M_n)$ be a weight sequence and suppose that
$\phi:[0,1] \rightarrow [0,1]$ induces a compact endomorphism $T$
of $D([0,1],M)$ with $\phi(x_0)=x_0.$
Then  $\sigma(T) \supset \{(\phi'(x_0))^n$: $n$ is a positive integer$\}\cup
\{0,1\}$.

{\bf Proof:} We first show that
$\phi'(x_0) \in \sigma(T).$ Indeed,
no $f \in D([0,1],M)$ satisfies $\phi'(x_0)f(x)-f(\phi(x))=x-x_0.$
For, if this were not the case and $f$ satisfied $\phi'(x_0)f(x)-
f(\phi(x))=x-x_0$, then $\phi'(x_0)f'(x)-f'(\phi(x))\phi'(x)=1.$
But evaluating at $x=x_0$ gives a contradiction. Hence $\phi'(x_0)
\in \sigma(T).$

Since $T$ is compact every non-zero
element in $\sigma(T)$ is an eigenvalue of $T$.
Clearly we may assume
that $\phi'(x_0)\neq 0$ and so, for some non-zero
$g \in D([0,1],M)$, $g(\phi(x))=\phi'(x_0)g(x).$ Therefore, for
each positive integer $n$, $g^n(\phi(x))=\phi'(x_0)^n g^n(x)$, whence
each $\phi'(x_0)^n \in \sigma(T)$ for all positive integers $n$.
\QED
     We remark that this lemma holds for non-compact endomorphisms
as well. For, by successive differentiation, one can show that if
$\phi(x_0)=x_0$, then for each positive integer $n$ there is no
$f \in D([0,1],M)$ satisfying $(\phi'(x_0))^n f(x) - f(\phi(x))=
(x-x_0)^n$. This again implies that $(\phi'(x_0))^n \in \sigma(T)$ for all
positive integers $n.$


{\bf Lemma 2.2:} Suppose that $(M_n)$ is a weight sequence,
$\phi:[0,1] \rightarrow [0,1]$
induces an endomorphism $T$ of $D([0,1],M)$,
$x_0$ is a fixed point of $\phi$ and $f$ is an eigenfunction of $T$
with eigenvalue $\lambda$.   If $\lambda \neq 0, 1$, and if
$\lambda \neq (\phi'(x_0))^n$
for all positive integers $n$, then $f^{(m)}(x_0)=0$ for all
non-negative integers $m$.

{\bf Proof:} Let $\lambda$ and $f$ satisfy the hypothesis.
Since $f$ is an eigenfunction,
 $\lambda f(x) = f(\phi(x))$,
and since $\lambda \neq 1,$ we have $f(x_0)=0.$

  Now suppose that $m \geq 1$ is given, and for $\nu < m$, $f^{(\nu)}(x_0)=0.$
We show that $f^{(m)}(x_0)=0.$ For each positive integer $m$,
$\lambda f^{(m)}(x)=f^{(m)}(\phi(x))(\phi'(x))^m+$[(terms containing sums
and products of $f^{(\nu)}$ and powers of $\phi^{(\nu)}$, $\nu < m.$)]
Since $f^{(\nu)}(x_0)=0$ for $\nu < m,$ evaluating at $x_0$, we have
$\lambda f^{(m)}(x_0)=f^{(m)}(x_0)(\phi'(x_0))^m,$ and since $\lambda
\neq (\phi'(x_0))^m,$ we conclude that $f^{(m)}(x_0)=0.$ Thus
$f^{(n)}(x_0)=0$ for all non-negative integers $n.$
\QED



{\bf Lemma 2.3:} Suppose that $\phi$ induces
a compact endomorphism of $D([0,1],M)$,
and that $x_0$ is the unique fixed point of $\phi$.
Let $\lambda$ be a non-zero complex number.
Suppose that $f \in D([0,1],M)$ is such that
$f(\phi(x))=\lambda f(x)$ and satisfies $f^{(k)}(x_0)=0$ for $k=0,1,2,\cdots.$
Then $f$ must be the zero function.
In particular, the endomorphism induced by $\phi$ has no non-zero eigenvalues other than
$\phi'(x_0)^n, n=0,1,2,\cdots$.

{\bf Proof:}
Suppose that $\lambda \neq 0$ and $f \in D([0,1],M)$ with $f(\phi(x))=\lambda f(x)$
and $f^{(k)}(x_0)=0$ for all $k=0,1,2,\dots .$ Since $|\phi'(x_0)| < 1,$
we can choose $\alpha$ with $|\phi'(x_0)| < \alpha < 1.$ Then there
exists $n_0 \in {\bf N}$ such that $\dis \phi_{n_0}$ (and hence all
later iterates of $\phi$) maps $[0,1]$ into a connected neighborhood $\cal U$
of $x_0$ with $|\phi'(t)|<\alpha$ for all $t \in \cal U.$
It then follows that there exists $C > 0$ such that
$\dis |\phi_n(x)-x_0| < \alpha^{n-n_0} < C \alpha^n$ for all $x \in (0,1).$
Next choose $m \in {\bf N}$ large enough that $\alpha^m < |\lambda|$. This
number $m$ will remain fixed for the remainder of this proof.
Applying Taylor's theorem to the real and imaginary parts of
$f$ we see that there is a constant $A>0$
such that $|f(x)| \leq A |x-x_0|^m$ for all $x \in [0,1]$.
Since $f \circ \phi = \lambda f,$ we have (for $\lambda \neq 0$),
that $f(x)=f(\phi_n(x))/\lambda^n$ for all
$x \in [0,1]$ and all $n \in {\bf N}$. Thus, for such $x$ and $n$, we have
$$|f(x)| = |f(\phi_n(x))/\lambda^n|
\leq A|\phi_n(x)-x_0|^m/|\lambda|^n
\leq AC(\alpha^n)^m/|\lambda|^n
=AC(\alpha^m/|\lambda|)^n.$$
Letting $n$ tend to infinity gives us $f(x)=0$, as required.
The last part of the conclusion now follows from Lemma 2.2.
\QED

{\bf Theorem 2.4:} Suppose that $(M_n)$ is a weight sequence,
$\phi$ induces a compact endomorphism $T$
of $D([0,1],M)$ and $x_0$ is the unique fixed point of $\phi$.
Then the spectrum
$\sigma(T)=
\{(\phi'(x_0))^n$: $n$ is a positive integer$\}\cup \{0,1\}$.

{\bf Proof:} Since $T$ is compact, $0 \in \sigma(T),$ and from
Theorem 1.1, $|\phi'(x_0)| < 1.$
Also
Lemma 2.1 shows that
$\sigma(T) \supset  \{(\phi'(x_0))^n$: $n$ is a positive integer$\}\cup
\{0,1\}$.
On the other hand since $T$ is compact
every non-zero element of $\sigma(T)$ is an eigenvalue. From
this and Lemma 2.3 it follows that
$\sigma(T) \subset \{(\phi'(x_0))^n$: $n$ is a positive integer$\}\cup
\{0,1\}$, and thus equality holds.
\QED
   We conclude this part by showing that the multiplicity of every
non-zero eigenvalue is $1.$

{\bf Theorem 2.5:} Suppose that $(M_n)$ is a weight sequence and that
$\phi$ induces a compact endomorphism $T$
of $D([0,1],M)$.
Then every non-zero eigenvalue of $T$ has multiplicity $1.$

{\bf Proof:} We start by showing that the eigenvalue $1$ has multiplicity
$1.$
Certainly $f(x) \equiv 1$ is an eigenvector. Next suppose
that $g \in D([0,1],M)$
is another  eigenvector. Then $g(\phi(x))=g(x)$
for all $x \in [0,1]$. It
follows easily that $g(\phi_n(x))=g(x)$ for all $x \in [0,1].$
However, since the induced endomorphism is compact, $\phi_n(x)$
converges to the unique fixed point, $x_0$, of $\phi$.
Hence $g(x) \equiv g(x_0)$
as claimed.

For the remainder of the proof we may assume that $\phi'(x_0) \neq 0$.
     Next suppose that $p$ is a positive integer. We show
that the eigenvalue $\phi'(x_0)^p$ also has multiplicity $1$.
We recall that  $|\phi'(x_0)| < 1.$

   We first show that if $f \in D([0,1],M)$ with
$f \circ \phi = (\phi'(x_0))^p f$, and if $f^{(p)}(x_0)=0$,
then $f \equiv 0.$ To this end assume that $f(\phi(x))=(\phi'(x_0))^p f(x)$
for all $x$. Then it is easy to see that $f^{(\nu)}(x_0)=0$
for $\nu=0,1, \cdots, p-1.$ Also when we assume that $f^{(p)}(x_0)=0,$
it is easy to see that $f^{(\nu)}(x_0) = 0$ for $\nu > p.$
Hence by Lemma 2.3, if $f \in D[0,1],M)$ and $f \circ \phi =
(\phi'(x_0))^p f$ and $f^{(p)}(x_0)=0$, then $f \equiv 0.$

   Finally suppose that $f, g \in D([0,1],M)$ and that both
$f$ and $g$ are non-zero eigenvectors for $\phi'(x_0)^p.$
For $f$ and $g$
to be non-zero, both $f^{(p)}(x_0) \neq 0$ and $g^{(p)}(x_0) \neq 0.$
Then the function $F$ defined by
$\dis F(x)=f(x)-\frac{f^{(p)}(x_0)}{g^{(p)}(x_0)}g(x)$ is an eigenvector
for $\phi'(x_0)^p$ and
satisfies $F^{(p)}(x_0)=0$.
Hence $F \equiv 0$ showing that $f$ and $g$ are
linearly dependent.
\QED

{\bf Part III. General endomorphisms}

     As noted in the introduction (Theorem A), 
the calculations from \cite{blaub} show that if
$(M_n)$ is a non-analytic sequence, then every analytic selfmap $\phi$ of $[0,1]$
with
$\|\phi'\|_\infty < 1$ induces a (compact) endomorphism 
of $D([0,1],M)$.
Also (Theorem A$'$), if
\[\frac{M_m}{m!}\frac{n!}{M_n}\leq \frac{B}{m^{n-m}},\]
for $m \leq n,$ and $\dis \frac{\|\phi^{(k)}\|_\infty}{k!}$ is
bounded, then
$\dis \|\phi'\|_\infty \leq 1$ is sufficient for $\phi$ to induce an endomorphism.

    We now relax the conditions on both $(M_n)$ and $\phi$ to prove a stronger version
of Theorem A$'$. The new condition on $(M_n)$ can be shown to be equivalent to
$\dis \sup_{n} \frac{n^2M_{n-1}}{M_n} < \infty,$ an easier condition to verify
than $\dis \frac{M_m}{m!}\frac{n!}{M_n}\leq \frac{B}{m^{n-m}}.$ In addition, $\phi$
can be assumed to be analytic, rather than satisfy the stronger condition
$\dis \sup_{k} \frac{\|\phi^{(k)}\|_\infty}{k!} < \infty.$

    We begin with two easily checked observations.

{\bf Observation 1:} If $(M_n)$ is a weight sequence, then
$\dis \sup_{n} \frac{n^2M_{n-1}}{M_n} < \infty$ if, and only if, there exist
positive constants $B$ and $D$ such that whenever $m \leq n$ we have
\[\dis \frac{M_m}{m!}\frac{n!}{M_n} \leq
B\left(\frac{D}{m}\right)^{n-m}.\] 

{\bf Observation 2:} If $(M_n)$ is a weight sequence and $a_j$ are non-negative
integers with $a_1+2a_2 + \cdots + na_n=n,$ then
\[\prod_{k=1}^n \left(\frac{M_k}{k!}\right)^{a_k}
\leq
\frac{M_n}{n!}.\]

{\bf Theorem 3.1:} Suppose that $(M_n)$ is a non-analytic weight sequence satisfying
$\dis \sup_{n} \frac{n^2M_{n-1}}{M_n} < \infty$ and that $\phi$ is a analytic
selfmap of $[0,1]$. If $\|\phi'\|_\infty \leq 1,$ then $\phi$ induces an endomorphism
of $D([0,1],M).$

{\bf Proof:} As noted, if $\dis \sup_{n} \frac{n^2M_{n-1}}{M_n} < \infty,$ then
there exist positive constants $B$ and $D$ such that $\dis \frac{M_m}{m!}\frac{n!}{M_n} \leq
B(\frac{D}{m})^{n-m}$, $m \leq n.$
 Suppose $\phi$ satisfies the hypotheses with $\|\phi'\|_{\infty} \leq 1.$  Let $F \in D([0,1],M).$ We show that $F \circ \phi \in D([0,1],M).$

The following equality for higher derivatives of composite functions is known
as Fa\`{a} di Bruno's formula.

\[\dis \frac{d^n}{dx^n}(F \circ \phi)=\sum_{m=0}^{n} F^{(m)}(\phi)
\left({\sum_{a}} \frac{n!}{a_1!a_2!\cdots a_n!} \prod_{k=1}^n{\left(\frac{\phi^{(k)}}{k!}\right)^{a_k}}\right)\]
where the inner sum $\dis\sum_a$ is over non-negative integers $a_1, a_2, \cdots,
a_n$ satisfying $a_1+a_2+ \cdots +a_n=m$ and $a_1+2a_2+\cdots +na_n=n.$

{\bf Throughout the proof  the inner sum $\dis\sum_{a}$ will always
be over non-negative integers $a_1, a_2, \cdots,a_n$ satisfying $a_1+a_2+ \cdots +a_n=m$ and $a_1+2a_2+\cdots +na_n=n.$}

    Thus, Fa\`{a} di Bruno's formula implies that

\[\dis \|\frac{d^n}{dx^n}(F \circ \phi) \|_{\infty} \leq \sum_{m=0}^n \|F^{(m)}
(\phi)\|_\infty \left({\sum_{a}}\frac{n!}{a_1!a_2!\cdots a_n!}
\prod_{k=1}^n{\left(\frac{\|\phi^{(k)}\|_\infty}{k!}\right)^{a_k}}\right)\]
and so

\[\dis \sum_{n=0}^{\infty} \frac{1}{M_n}\|\frac{d^n}{dx^n}(F \circ \phi)\|
_{\infty} \leq \sum_{n=0}^{\infty} \frac{1}{M_n} \sum_{m=0}^n \|F^{(m)}\|
_{\infty} \left({\sum_{a}}\frac{n!}{a_1!a_2!\cdots a_n!}
\prod_{k=1}^n{\left(\frac{\|\phi^{(k)}\|_\infty}{k!}\right)^{a_k}}\right).\]

Then, after interchanging the order of summation, we have

\[ \sum_{n=0}^{\infty} \frac{1}{M_n}\|\frac{d^n}{dx^n}(F \circ \phi)\|
_{\infty} \leq \sum_{m=0}^{\infty} \|F^{(m)}\|_\infty \sum_{n=m}^{\infty} \frac{1}{M_n}
\left({\sum_{a}}\frac{n!}{a_1!a_2!\cdots a_n!}
\prod_{k=1}^n{\left(\frac{\|\phi^{(k)}\|_\infty}{k!}\right)^{a_k}}\right).\]
Since $\phi$ is analytic, there is some $C>0$ such that
$\dis \frac{\|\phi^{(k)}\|_\infty}{k!} \leq C^k$ for
all $k$. We then have (denoting the greatest integer less than or
equal to $x$ by $[x]$) that
$\dis \sum_{n=0}^\infty \frac{1}{M_n}\|(F \circ \phi)^{(n)}\|_\infty \leq
\dis {\cal A}_1 + {\cal A}_2$ where

\[\dis {\cal A}_1 = \sum_{m=0}^{[CD]}\|F^{(m)}\|_\infty
\sum_{n=m}^{\infty} \frac{1}{M_n}
\left({\sum_{a}}\frac{n!}{a_1!a_2!\cdots a_n!}
\prod_{k=1}^n{\left(\frac{\|\phi^{(k)}\|_\infty}{k!}\right)^{a_k}}\right)\]
and
\[\dis {\cal A}_2 =\sum_{m=[CD]+1}^{\infty}\|F^{(m)}\|_\infty
\sum_{n=m}^{\infty} \frac{1}{M_n}
\left({\sum_{a}}\frac{n!}{a_1!a_2!\cdots a_n!}
\prod_{k=1}^n{\left(\frac{\|\phi^{(k)}\|_\infty}{k!}\right)^{a_k}}\right).\]

     First we estimate ${\cal A}_1.$
Recalling that, for the non-negative integers $a_k$ under consideration,
\[\prod_{k=1}^n \left(\frac{M_k}{k!}\right)^{a_k}
\leq
\frac{M_n}{n!},\]
we see that, for each $m$, the coefficient of
$\|F^{(m)}\|_\infty$ is

\[\sum_{n=m}^{\infty} \frac{1}{M_n}
\left({\sum_{a}}\frac{n!}{a_1!a_2!\cdots a_n!}
\prod_{k=1}^n{\left(\frac{\|\phi^{(k)}\|_\infty}{k!}\right)^{a_k}}\right)\]


\[=\sum_{n=m}^\infty \frac{1}{n!} \left({\sum_{a}}\frac{n!}{a_1!a_2!\cdots a_n!}
\frac{n!}{M_n} \prod_{k=1}^n{\left(\frac{M_k}{k!}
\frac{\|\phi^{(k)}\|_\infty}{M_k}\right)^{a_k}}\right)\]

\[\leq \sum_{n=m}^\infty \frac{1}{n!}
\left({\sum_{a}}\frac{n!}{a_1!a_2!\cdots a_n!}
\prod_{k=1}^n{\left(\frac{\|\phi^{(k)}\|_\infty}{M_k}\right)^{a_k}}\right).\]

Now from formula B3, page 823 in \cite{abr}, the right hand side of the last
inequality equals $\dis \frac{1}{m!}
\left(\sum_{k=1}^\infty \frac{\|\phi^{(k)}\|_\infty}
{M_k}\right)^m=\frac{1}{m!}\left(\|\phi\|_{D}-\|\phi\|_\infty\right)^m.$

Therefore, 
\[\dis {\cal A}_1=
\sum_{m=0}^{[CD]}\|F^{(m)}\|_\infty
\sum_{n=m}^{\infty} \frac{1}{M_n}
\left({\sum_{a}}\frac{n!}{a_1!a_2!\cdots a_n!}
\prod_{k=1}^n{\left(\frac{\|\phi^{(k)}\|_\infty}{k!}\right)^{a_k}}\right)\]
\[\dis \leq \sum_{m=0}^{[CD]} 
\frac{\|F^{(m)}\|_\infty}{m!}\left(\|\phi\|_{D}-\|\phi\|_\infty\right)^m,\]
which is finite.

     We next show that ${\cal A}_2$ is finite. Given that
$\dis \frac{n!}{M_n}\frac{M_m}{m!} \leq B\left(\frac{D}{m}\right)^{n-m},$
we have
\[\dis {\cal A}_2=\sum_{m=[CD]+1}^{\infty}\|F^{(m)}\|_\infty
\sum_{n=m}^{\infty} \frac{1}{M_n}
\left({\sum_{a}}\frac{n!}{a_1!a_2!\cdots a_n!}
\prod_{k=1}^n{\left(\frac{\|\phi^{(k)}\|_\infty}{k!}\right)^{a_k}}\right)\]
\[\dis=\sum_{[CD]+1}^{\infty}\frac{\|F^{(m)}\|_\infty}{M_m}m!
\sum_{n=m}^{\infty} \frac{1}{n!}
\left({\sum_{a}}\frac{n!}{M_n}\frac{M_m}{m!}\frac{n!}{a_1!a_2!\cdots a_n!}
\prod_{k=1}^n{\left(\frac{\|\phi^{(k)}\|_\infty}{k!}\right)^{a_k}}\right)\]
\[\dis \leq B \sum_{[CD]+1}^{\infty}\frac{\|F^{(m)}\|_\infty}{M_m} m!
\sum_{n=m}^{\infty} \frac{1}{n!}
\left({\sum_{a}}\frac{n!}{a_1!a_2!\cdots a_n!}
\prod_{k=1}^n{\left(\frac{\|\phi^{(k)}\|_\infty D^{k-1}}{k! m^{k-1}}\right)^{a_k}}\right)\]
\[=\dis B \sum_{[CD]+1}^\infty \frac{\|F^{(m)}\|_\infty}{M_m}\left(\sum_{k=1}^\infty
\frac{\|\phi^{(k)}\|_\infty D^{k-1}}{k!m^{k-1}}\right)^m,\]
where we have again used
the previously quoted formula from \cite{abr}.
Therefore (assuming, as we may, that $\phi$ is non-constant),

\[\dis {\cal A}_2 \leq B \sum_{[CD]+1}^\infty
\frac{\|F^{(m)}\|_\infty}{M_m}\|\phi'\|_\infty ^m
\left(1 + \sum_{k=2}^\infty \frac{\|\phi^{(k)}\|_\infty D^{k-1}}{k!\|\phi'\|_\infty m^{k-1}}\right)^m\]
\[=\dis B \sum_{[CD]+1}^\infty \frac{\|F^{(m)}\|_\infty}{M_m}\|\phi'\|_\infty^m
\left(1 + \sum_{k=2}^\infty \frac{C^kD^{k-1}}{\|\phi'\|_\infty m^{k-1}}\right)^m\]
\[\dis =  B \sum_{[CD]+1}^\infty \frac{\|F^{(m)}\|_\infty}{M_m}\|\phi'\|_\infty^m
\left(1 + \frac{C}{\|\phi'\|_\infty}\sum_{k=2}^\infty \frac{C^{k-1} D^{k-1}}{m^{k-1}}\right)^m\]
\[\dis =  B \sum_{[CD]+1}^\infty \frac{\|F^{(m)}\|_\infty}{M_m}\|\phi'\|_\infty^m
\left(1 + \frac{C^2D}{\|\phi'\|_\infty}\frac{1}{m-CD}\right)^m.\]

    Since $\dis \left(1 + \frac{C^2D}{\|\phi'\|_\infty}\frac{1}{m-CD}\right)^m$ is easily seen
to be bounded in $m$, we conclude that ${\cal A}_2$ is finite when
$\|\phi'\|_{\infty} \leq 1,$ as claimed.
\QED

Most non-analytic weight sequences whose growth is, in some sense, regular and which do
not satisfy the conditions of Theorem 3.1 are such that
$n^2M_n/M_{n+1} \to \infty$ as $n \to \infty$.
In this case $\|\phi'\|_\infty \leq 1$ is no longer sufficient for $\phi$
to induce an endomorphism of $D([0,1],M)$ as the following theorem shows.

{\bf Theorem 3.2:}
Let $(M_n)$ be a weight sequence such that
$n^2M_n/M_{n+1} \to \infty$ as $n \to \infty$.
Then the map $\phi(x)=(1+x^2)/2$ does not induce an endomorphism
of $D([0,1],M)$.

{\bf Proof:}
Consider the functions $F_A (t) = \exp(A(t-1))$ where $A$ varies through the
positive real numbers. Clearly when $t=1$, all the derivatives
involved
are non-negative.

Then for $n \geq 2$, the second
term $(m=n-1)$ of the Fa\`a di Bruno expansion of $(F_A \circ \phi)^{(n)}(1)$
tells us that
$$(F_A \circ \phi)^{(n)}(1)
\geq {n \choose 2} F_A^{(n-1)}(\phi(1)) (\phi'(1))^{n-2} \phi''(1)$$
$$= {n \choose 2} F_A^{(n-1)}(1) = {n \choose 2} A^{n-1}.$$

Set $C=1/4$.
We obtain, for $n \geq 2$, and such functions $F_A$,
$\norm{(F_A \circ \phi)^{(n)}}_\infty \geq C n^2 \norm{F_A^{(n-1)}}_\infty$.

Let $B>0$.
Choose $N \in {\bf N}$ with $N > 2$ and such that, for all $k \geq N$ we
have
$k^2 M_{k-1}/M_k > B$.
Choose $A > 0$ such that
$\sum_{k=0}^{N-1}{ A^k/M_k < (1/2)A^N/M_N}$.
Then, for the function $F_A$ we
have
$$\sum_{k=N}^\infty{\norm{F_A^{(k)}}_\infty/M_k} =
\sum_{k=N}^\infty{A^k/M_k}$$
$$\geq (1/2) \norm{F_A}_{D}.$$
But, by the above, we also have
$$\norm{F_A \circ \phi}_{D}
  \geq \sum_{k=N}^\infty \norm{(F_A\circ \phi)^{(k)}}_\infty/M_k$$
$$\geq \sum_{k=N}^\infty {C{k^2} \norm{F_A^{(k-1)}}_\infty/M_k}
=C \sum_{k=N}^\infty {{k^2}  A^{k-1}/M_k}$$
$$=C \sum_{k=N}^\infty {{k^2} (M_{k-1}/M_k) (A^{k-1}/M_{k-1})}
\geq C B \sum_{k=N}^\infty { A^{k-1}/M_{k-1}}$$
$$\geq (C/2) B \norm{F_A}_{D}.$$

As $B>0$ was arbitrary, this shows that  $\phi$
cannot induce a bounded endomorphism of $D([0,1],M)$ and so $\phi$
does not induce an endomorphism at all.
\QED

   The rate of growth of $\dis \frac{n^2M_n}{M_{n+1}}$ and the order
of the entire function $\dis g(z)=\sum^\infty_{n=0}\frac{z^n}{M_n}$ are
somewhat related. In particular, if $\dis \frac{n^2M_n}{M_{n+1}}$ is bounded,
then the order $\rho$ of $g$ is less than $\frac{1}{2},$ or $\rho=\frac{1}{2}$
and the type $\tau$ is finite. On the other hand, if
$\dis \lim_{n \rightarrow \infty} \frac{n^2M_n}{M_{n+1}}=\infty$, then
$\rho > \frac{1}{2}$, or $\rho = \frac{1}{2}$ and $\tau = \infty.$
Using order and type of the function $g(z)$, one can also show the following.

{\bf Proposition 3.3:} Let $(M_n)$ be a non-analytic weight function.
Suppose that $\phi:[0,1] \rightarrow [0,1]$ is in $D([0,1],M)$
with $\phi(1)=1=\phi'(1)=\|\phi'\|_\infty$,  and
$\phi^{(\nu)}(1) \geq 0$ for all $\nu.$
Let $\dis g(z)=\sum_{n=0}^\infty \frac{z^n}{M_n},$ an entire function
with order $\rho$ and type $\tau.$
If $\phi''(1) \neq 0$ and $\rho > \frac{1}{2},$ then
$\phi$ does not induce an endomorphism of $D([0,1],M).$
Also, if $ \rho=\frac{1}{2}$ and $\tau=\infty$, then such
$\phi$ does not induce an endomorphism.



\vspace{.3in}

{\sf  School of Mathematical Sciences

 University of Nottingham

 Nottingham NG7 2RD, England

 email: Joel.Feinstein@nottingham.ac.uk

and

 Department of Mathematics

 University of Massachusetts at Boston

 100 Morrissey Boulevard

 Boston, MA 02125-3393

 email: hkamo@cs.umb.edu

\vspace{.3in}
2000 Mathematics Subject Classification: 46J15, 47B48
\vspace{.3in}

This research was supported by EPSRC grants GR/M31132 and GR/R09589 }

\end{document}